\documentclass{tran-l}

 \newtheorem{thm}{Theorem}[subsection]
 \newtheorem{cor}[thm]{Corollary}
 \newtheorem{lem}[thm]{Lemma}
 \newtheorem{prop}[thm]{Proposition}
 \theoremstyle{definition}
 \newtheorem{defn}[thm]{Definition}
 \theoremstyle{remark}
 \newtheorem{rem}[thm]{Remark}
 \numberwithin{equation}{subsection}
\newcommand{\pn}{\noindent}

\newcommand{\ZZ}{\mathbb{Z}}

\newcommand{\CC}{\mathbb{C}}
\newcommand{\GG}{\mathbb{G}}

\newcommand{\Hom}{\mathrm{Hom}}

\newcommand{\Biext}{\mathrm{Biext}}

\newcommand{\bBiext}{\mathrm{\mathbf{Biext}}}
\newcommand{\W}{\mathrm{W}}

\newcommand{\Gr}{\mathrm{Gr}}
\newcommand{\rk}{\mathrm{rk}}

\begin{document}

\title[Biextensions and 1-motives]
{Biextensions of 1-motives by 1-motives}

\author{Cristiana Bertolin}

\address{NWF I -Mathematik, Universit\"at Regensburg, 93040 Regensburg,
 Germany}

\email{cristiana.bertolin@mathematik.uni-regensburg.de}



\date{}
\dedicatory{}

\commby{Cristiana Bertolin}

%%% ----------------------------------------------------------------------

\begin{abstract}
Let $S$ be a scheme. In this paper, we define the notion of biextensions of 1-motives by 1-motives.
Moreover, if ${\mathcal{M}}(S)$ denotes the Tannakian category generated by 1-motives over $S$ (in a geometrical sense), we define geometrically the morphisms of ${\mathcal{M}}(S)$
from the tensor product of two 1-motives $M_1 \otimes M_2$ to another 1-motive $M_3,$ to be the isomorphism classes of biextensions of $(M_1,M_2)$ by $M_3$:

\[ {\Hom}_{{\mathcal{M}}(S)}(M_1 \otimes M_2,M_3) = {\Biext}^1(M_1,M_2;M_3) \]

\pn Generalizing this definition we obtain, modulo isogeny, the geometrical notion of morphism of ${\mathcal{M}}(S)$ from a finite tensor product of 1-motives to another 1-motive.
\end{abstract}

%%% ----------------------------------------------------------------------
\maketitle
%%% ----------------------------------------------------------------------

\section*{Introduction}
Let $S$ be a scheme. 

A 1-motive over $S$ consists of a $S$-group scheme $X$ which is locally for the \'etale topology a constant group scheme defined by a finitely generated free
$\ZZ$-module, a semi-abelian $S$-scheme $G$, and a morphism $u:X \longrightarrow G$ of $S$-group schemes.

Let ${\mathcal{M}}(S)$ be what should be the Tannakian category generated by 1-motives over $S$ in a geometrical sense. We know very little about this category ${\mathcal{M}}(S)$: in particular, we are not able to describe geometrically the object of ${\mathcal{M}}(S)$ defined as the tensor product of two 1-motives! 
Only if $S={\mathrm{Spec}}\,(k),$ with $k$  a field of characteristic 0 embeddable in $\CC$, we know something about ${\mathcal{M}}(k)$:
in fact, identifying 1-motives with their mixed realizations, we can identify ${\mathcal{M}}(k)$ with the Tannakian sub-category of an ``appropriate'' Tannakian category of mixed realizations generated by the mixed realizations of 1-motives.

The aim of this paper is to use biextensions in order to define some morphisms 
in the category ${\mathcal{M}}(S)$: \emph{Geometrically the morphisms of ${\mathcal{M}}(S)$ from the tensor product of two 1-motives $M_1 \otimes M_3$ to another 1-motive $M_3$ are the isomorphism classes of biextensions of $(M_1,M_2)$ by $M_3$}:

\begin{equation}\label{defmorintro}
{\Hom}_{{\mathcal{M}}(S)}(M_1 \otimes M_2,M_3) = {\Biext}^1(M_1,M_2;M_3) 
\end{equation}

The idea of defining morphisms through biextensions goes back to Alexander Gro\-then\-dieck. In fact, in~\cite{SGA7} Expos\'e VIII he defines pairings from biextensions:
if $P,Q,G$ are three abelian groups of a topos \textbf{T}, to each isomorphism class of biextensions of $(P,Q)$ by $G,$ he associates a pairing $({}_{l^n}P)_{n \geq 0} \otimes ({}_{l^n}Q)_{n \geq 0}  \longrightarrow ({}_{l^n}G)_{n \geq 0} $ where $({}_{l^n}P)_{n \geq 0} $ (resp. $({}_{l^n}Q)_{n \geq 0} $, $({}_{l^n}G)_{n \geq 0} $) is the projective system constructed from the kernels ${}_{l^n}P$ ( resp. 
${}_{l^n}Q$ , ${}_{l^n}G$) of the multiplication by $l^n$ for each $n \geq 0$.

Generalizing Grothendieck's work, in ``Theorie de Hodge III'' Pierre Deligne defines the notion of biextension of two complexes of abelian groups concentrated in degree 0 and -1 (over any topos \textbf{T}) by an abelian group. Applying this definition to 1-motives and to ${\GG}_m,$ to each isomorphism class of such biextensions he associates a pairing ${\mathrm{T}}_*(M_1) \otimes {\mathrm{T}}_*(M_2)
\longrightarrow {\mathrm{T}}_*( {\GG}_m)$ in the Hodge, De Rhan, $\ell$-adic realizations (resp. $*={\mathrm{H}}, *={\mathrm{dR}}, *=\ell$).

Our definition of biextension of 1-motives by 1-motive generalizes the one of Deligne. The main idea of our definition is the following: Recall that a 1-motive $M$ over $S$ can be described also as a 7-uplet $(X,Y^{\vee},A,A^*, v ,v^*,\psi)$ where

\begin{itemize}
    \item $X$ and $Y^{\vee}$ are two $S$-group schemes which are locally 
for the \'etale topology constant group schemes defined by 
 finitely generated free $\ZZ$-modules;
    \item $A$ and $A^*$ are two abelian $S$-schemes dual to each other;
    \item $v:X \longrightarrow A$ and $v^*:Y^{\vee} \longrightarrow A^*$ are two morphisms  of $S$-group schemes; and 
    \item $\psi$ is a trivialization of the pull-back $(v,v^*)^*{\mathcal{P}}_A$ via $(v,v^*)$ of the Poincar\' e biextension ${\mathcal{P}}_A$ of $(A,A^*)$.
\end{itemize}    

\pn In other words, we can reconstruct the 1-motive $M$ from the Poincar\' e biextension ${\mathcal{P}}_A$ of $(A,A^*)$ and some trivializations of some pull-backs of ${\mathcal{P}}_A$. The 4-uplet $(X,Y^{\vee},A,A^*)$ corresponds to the pure part of the 1-motive, i.e. it defines the pure motives underlying $M$, and the 3-uplet
$v ,v^*,\psi$ represents the ``mixity'' of $M$.  
Therefore the Poincar\' e biextension ${\mathcal{P}}_A$ is related to the pure part 
$(X,Y^{\vee},A,A^*)$ of the 1-motive and the trivializations are related to the mixed part $v ,v^*,\psi$ of $M$.\\
Now let $M_i=[X_i \stackrel{u_i}{\longrightarrow}G_i]$ (for $i=1,2,3$) be a 1-motive over $S$. A biextension of $(M_1,M_2)$ by $M_3$ is a biextension $\mathcal{B}$ of   
$(G_1,G_2)$ by $G_3$, some trivializations of some pull-backs of $\mathcal{B}$, and a morphism $X_1 \otimes X_2 \longrightarrow X_3$ compatible with the above trivializations.
Since by the main Theorem of~\cite{B2}, to have a biextension of $(G_1,G_2)$ by $G_3$
is the same thing has to have a biextension of $(A_1,A_2)$ by $Y_3(1)$, we can simplify our definition: a biextension $B=(B,\Psi_1, \Psi_2,\Psi,\lambda)$ of $(M_1,M_2)$ by $M_3$ consists of

\begin{enumerate} 
  \item a biextension of $B$ of $(A_1, A_2)$ by $Y_3(1);$
   \item a trivialization $\Psi_1 $ (resp. $\Psi_2$) of the biextension 
$(v_1, id_{A_2})^*B$ (resp. \\
 $(id_{A_1},v_2)^*B$) of $(X_1, A_2)$ by $Y_3(1)$ 
(resp. of $(A_1, X_2)$ by $Y_3(1)$) obtained as pull-back of the biextension $B$ via $(v_1,id_{A_2})$ (resp. via $(id_{A_1},v_2)$); 
   \item a trivialization $\Psi$ of the biextension 
$(v_1, v_2)^*B$ of $(X_1, X_2)$ by $Y_3(1)$ obtained as pull-back of the biextension 
$B$ via $(v_1,v_2)$, which coincides with the trivializations 
induced by $\Psi_1$ and $\Psi_2$ over $X_1 \times X_2$;
    \item a morphism $\lambda:X_1 \otimes X_2 \longrightarrow X_3$ compatible with the trivialization $\Psi$ and the morphism $u_3:X_3 \longrightarrow G_3$. 
\end{enumerate}    

\pn The biextension $B$ and the morphism $\lambda:X_1 \otimes X_2 \longrightarrow X_3$ take care of the pure parts of the 1-motives $M_1,M_2$ and $M_3$, i.e. of
$X_i,Y_i,A_i,A^*_i$ for $i=1,2,3$. On the other hand the trivializations $\Psi_1, \Psi_2$, $\Psi,$ and the compatibility of $\lambda$ with such trivialirations and with $u_3,$ take care of the mixed parts of $M_1,M_2$ and $M_3$, i.e. of $v_i,v^*_i$ and $\psi_i$ for $i=1,2,3$. \\
In terms of morphisms, the biextension $B=(B,\Psi_1, \Psi_2,\Psi,\lambda)$ defines a morphism $M_1 \otimes M_2 \longrightarrow M_3$  where

\begin{itemize}
    \item the biextension $B$ of $(A_1,A_2)$ by $Y_3(1)$ defines the component $A_1 \otimes A_2 \longrightarrow Y_3(1)$;
    \item the morphism $\lambda$ defines the component $X_1 \otimes X_2 \longrightarrow X_3$ of weight 0;
    \item the trivializations $\Psi_1, \Psi_2$, $\Psi$ and the compatibility of $\lambda$ with such trivializations and with $u_3$ take care 
of the other components of $M_1 \otimes M_2 \longrightarrow M_3$, and
of the compatibility of the mixed part of $M_1 \otimes M_2$ and the mixed part of $M_3$ through this morphism of $M_1 \otimes M_2 \longrightarrow M_3$.
\end{itemize}

At the end of Section 1, we give several examples of biextensions of 1-motives by 1-motives. At the end of Section 2, we show how to construct explicitly the morphisms corresponding to the biextensions stated as example at the end of Section 1.

Then we verify that the property of \emph{respecting weights} which is satisfied by the morphisms of ${\mathcal{M}}(S)$, is also verified by biextensions.
 For example, in ${\mathcal{M}}(S)$ there are no morphisms from $X_1 \otimes X_2 $ to $G_3$ or from $G_1 \otimes G_2$ to $X_3$
(Lemma \ref{georem0}), and therefore we must have that all the biextensions of $([X_1\rightarrow 0],[X_2\rightarrow 0])$ by $[0 \rightarrow G_3]$ and all the biextensions $([0 \rightarrow G_1],[0 \rightarrow G_2])$ by $[X_3\rightarrow 0]$ are trivial. To have a morphism from $G_1 \otimes G_2$ to $G_3$ is the same thing as to have a morphism from $A_1 \otimes A_2$ to $Y_3(1)$ (Corollary \ref{georem2})
and therefore we must have that to have a biextension of  $(G_1, G_2)$ by $G_3$ is the same thing as to have a biextension of $(A_1, A_2)$ by $Y_3(1)$ (cf. Theorem \cite{B2}).

We can extend definition (\ref{defmorintro}) to a finite tensor product of 1-motives in the following way: again because of weights, a morphism from a finite tensor product $\otimes^l_1 M_j$ of 1-motives to a 1-motive $M$ involves only the quotient 
$\otimes^l_1 M_j / {\W}_{-3}(\otimes^l_1 M_j)$ of the mixed motive $\otimes^l_1 M_j$. However the motive $\otimes^l_1 M_j / {\W}_{-3}(\otimes^l_1 M_j)$ is isogeneous to
a finite sum of copies of $M_{\iota_1} \otimes M_{\iota_2}$ for $\iota_1, \iota_2 \in \{1, \dots,l\}$ (Lemma \ref{otimesM}), and therefore, modulo isogeny, a morphism 
from a finite tensor product of 1-motives to a 1-motive is a sum of isomorphism classes of biextensions of 1-motives by 1-motives (Theorem \ref{thmotimes}).

A special case of definition (\ref{defmorintro}) was already used in the computation of the unipotent radical of the Lie algebra of the motivic Galois group of a 1-motive defined over a field $k$ of characteristic 0 (cf.~\cite{B}).
In fact in~\cite{B} (1.3.1), using Deligne's definition of biextension 
of 1-motives by ${\GG}_m$, we define a morphism from the tensor product $M_1 \otimes M_2$ of two 1-motives  to a torus as an isomorphism class of biextensions of 
$(M_1,M_2)$ by this torus.

%------------------------------------------------------------------------

\section*{Acknowledgment}

I want to thank the referee of my article~\cite{B} who suggested me to generalize it: the study of biextensions of 1-motives by 1-motives started with the attempt of generalizing the definition~\cite{B} (1.3.1).
\vspace{0.5cm}

In this paper $S$ is a scheme.

%-------------------------------------------------------------------------

\section{Biextensions of 1-motives by 1-motives}

In~\cite{D1} (10.1.10), Deligne defines a \textbf{ 1-motive} $M$ over $S$ as

\begin{enumerate}
    \item  a $S$-group scheme $X$ which is locally for the \'etale
topology a constant group scheme defined by a finitely generated free
$\ZZ$-module,
    \item a semi-abelian $S$-scheme $G$, i.e. 
an extension of an abelian $S$-scheme $A$ by a $S$-torus $Y(1),$
with cocharacter group $Y$,
    \item a morphism $u:X \longrightarrow G$ of $S$-group schemes.
\end{enumerate}

The 1-motive $M$ can be view also as a complex $[X \stackrel{u}{\longrightarrow}G]$ of commutative $S$-group schemes concentrated in degree 0 and -1. 
An \textbf{isogeny between two 1-motives}
$M_{1}=[X_{1} {\buildrel u_{1} \over \longrightarrow} G_{1}]$ and 
$M_{2}=[X_{2} {\buildrel u_{2} \over \longrightarrow} G_{2}]$ is a morphism of complexes $(f_{X},f_{G})$ such that 
 $f_{X}:X_{1} \longrightarrow X_{2}$ is injective with finite cokernel, and 
 $f_{G}:G_{1} \longrightarrow G_{2}$ is surjective with finite kernel.

1-motives are mixed motives of niveau $\leq 1$: the weight filtration $W_*$ on $M=[X \stackrel{u}{\longrightarrow} G] $ is 

\begin{eqnarray}
\nonumber  {\W}_{i}(M) &=& M  ~~{\rm  for ~ each~} i \geq 0, \\
\nonumber  {\W}_{-1}(M) &=& [0 \longrightarrow G], \\
\nonumber  {\W}_{-2}(M) &=& [0 \longrightarrow  Y(1)], \\
\nonumber  {\W}_{j}(M) &=& 0 ~~~{\rm  for ~ each~} j \leq -3.
\end{eqnarray}

\pn In particular, we have  ${\rm Gr}_{0}^{{\W}}(M)= 
[X {\buildrel  \over \longrightarrow} 0], {\rm Gr}_{-1}^{{\W}}(M)= 
[0 {\buildrel  \over \longrightarrow} A]$ and $ {\rm Gr}_{-2}^{{\W}}(M)= 
[0 {\buildrel  \over \longrightarrow}  Y(1)].$

\subsection{The category of biextensions of 1-motives by 1-motives}
Let $M_i=[X_i \stackrel{u_i}{\longrightarrow}G_i]$ (for $i=1,2,3$) be a 1-motive over $S$. The following definition 
of biextension of $(M_1,M_2)$ by $M_3$, is a generalization of Deligne's 
definition~\cite{D1} (10.2):

\begin{defn}\label{defbiext}
A \textbf{biextension $\mathcal{B}= (\mathcal{B}, \Psi_1, \Psi_2,\Psi,\lambda)$ of $(M_1,M_2)$ by $M_3$} consists of

\begin{enumerate}
    \item a biextension of $\mathcal{B}$ of $(G_1, G_2)$ by $G_3$;
    \item a trivialization (= biaddictive section) $\Psi_1 $ (resp. $\Psi_2$) of the biextension $(u_1, id_{G_2})^* {\mathcal{B}}$ (resp. $(id_{G_1},u_2)^*{\mathcal{B}}$) of $(X_1, G_2)$ by $G_3$ (resp. $(G_1, X_2)$ 
by $G_3$) obtained as pull-back of the biextension ${\mathcal{B}}$ via $(u_1,id_{G_2})$ (resp. $(id_{G_1},u_2)$); 
    \item a trivialization $\Psi$ of the biextension $(u_1, u_2)^*{\mathcal{B}}$
 of $(X_1, X_2)$ 
by $G_3$ obtained as pull-back of the biextension 
${\mathcal{B}}$ by $(u_1,u_2)$, which coincides with the trivializations 
induced by $\Psi_1$ and $\Psi_2$
 over $X_1 \times X_2$, i.e. 
\[(u_1,id_{G_2})^*\Psi_2=\Psi=(id_{G_1},u_2)^*\Psi_1; \]
    \item a morphism $\lambda:X_1 \times X_2 \longrightarrow X_3$ of $S$-group schemes such that $u_3 \circ \lambda: X_1 \times X_2 \longrightarrow G_3$ is compatible with the trivialization  $\Psi$ of the biextension  $(u_1,u_2)^* {\mathcal{B}}$ of $(X_1, X_2)$ by $G_3$. 
\end{enumerate}    
\end{defn}

Let $M_i=[X_i \stackrel{u_i}{\longrightarrow}G_i]$ and 
$M'_i=[X'_i \stackrel{u'_i}{\longrightarrow}G'_i]$ (for $i=1,2,3$)
be 1-motives over $S$. Moreover let $({\mathcal{B}},\Psi_{1}, \Psi_{2},\lambda)$ be a biextension of $(M_1,M_2)$ by $M_3$ and let $({\mathcal{B}}',\Psi'_{1}, \Psi'_{2},\lambda')$ be a biextension of $(M'_1,M'_2)$ by $M'_3$.

\begin{defn}\label{defmorphbiext}
A \textbf{morphism of biextensions}
\[(F,\Upsilon_1,\Upsilon_2,\Upsilon,g_3):({\mathcal{B}},\Psi_{1}, \Psi_{2},\lambda)
\longrightarrow ({\mathcal{B}}',\Psi'_{1}, \Psi'_{2},\lambda')\]
\pn consists of
\begin{enumerate}
    \item a morphism $F=(F,f_1,f_2,f_3):{\mathcal{B}} \longrightarrow {\mathcal{B}}'$ from the biextension ${\mathcal{B}}$ to the biextension ${\mathcal{B}}'$. In particular,
\[ f_1:G_1 \longrightarrow G'_1 \qquad f_2:G_2 \longrightarrow G'_3 \qquad f_3:G_3 \longrightarrow G'_3\]
\pn are morphisms of groups $S$-schemes.
    \item a morphism of biextensions 

\[\Upsilon_1=(\Upsilon_1,g_1,f_2,f_3):(u_1, id_{G_2})^* {\mathcal{B}} \longrightarrow (u'_1, id_{G'_2})^* {\mathcal{B}}'\]

\pn compatible with the morphism $F=(F,f_1,f_2,f_3)$ and with the trivializations $\Psi_{1}$ and $\Psi'_{1}$, and  a morphism of biextensions 

\[\Upsilon_2=(\Upsilon_2,f_1,g_2,f_3):(id_{G_1},u_2)^* {\mathcal{B}} \longrightarrow (id_{G'_1},u'_2)^* {\mathcal{B}}'\]

\pn compatible with the morphism $F=(F,f_1,f_2,f_3)$ and with the trivializations $\Psi_{2}$ and $\Psi'_{2}$. In particular 
\[g_1: X_1 \longrightarrow X'_1 \qquad g_2: X_2 \longrightarrow X'_2  \]
\pn are morphisms of groups $S$-schemes.
    \item a morphism of biextensions $\Upsilon=(\Upsilon,g_1,g_2,f_3):
(u_1,u_2)^* {\mathcal{B}} \longrightarrow (u'_1, u'_2)^* {\mathcal{B}}'$ compatible with the morphism $F=(F,f_1,f_2,f_3)$ and with the trivializations $\Psi$ and $\Psi'$.
    \item a morphism $g_3:X_3 \longrightarrow X'_3$ of $S$-group schemes compatible with $u_3$ and $u'_3$ (i.e. $u'_3 \circ g_3 =f_3 \circ u_3$) and such that
\[ \lambda' \circ (g_1 \times g_2)= g_3 \circ \lambda\]
\end{enumerate}    
\end{defn}

\begin{rem}\label{remmorphbiext}
 The pair $(g_3,f_3)$ defines a morphism of ${\mathcal{M}}(S)$ from  $M_3$ to $M'_3$. The pairs $(g_1,f_1)$ and $(g_2,f_2)$ define morphisms of ${\mathcal{M}}(S)$ from $M_1$ and $M'_1$ and from $M_2$ to $M'_2$ respectively.
\end{rem}

We denote by ${\bBiext}(M_1,M_2;M_3)$ the category of biextensions 
of $(M_1,M_2)$ by $M_3$. Like in the category of biextensions of semi-abelian schemes,
in the category ${\bBiext}(M_1,M_2;M_3)$ we can make the sum of two objects.
Let ${\Biext}^0(M_1,M_2;M_3)$ be the group of automorphisms of any biextension 
of $(M_1,M_2)$ by $M_3$, and let\\
 ${\Biext}^1(M_1,M_2;M_3)$ be the group of isomorphism classes of biextensions of $(M_1,M_2)$ by $M_3$.

\subsection{A more useful definition}
From now on we will work on the topos ${\mathrm{\mathbf{T}_{fppf}}}$ 
associated to the site of locally of finite presentation $S$-schemes,
 endowed with the fppf topology.

Proposition~\cite{D1} (10.2.14) furnishes a more symmetric description of 1-motives: 
consider the 7-uplet $(X,Y^{\vee},A,A^*, v ,v^*,\psi)$ where

\begin{itemize}
    \item $X$ and $Y^{\vee}$ are two $S$-group schemes which are locally 
for the \'etale topology constant group schemes defined by 
 finitely generated free $\ZZ$-modules. We have to think at $X$ and at $Y^{\vee}$ as character groups of $S$-tori which should be written (according to our notation) $X^{\vee}(1)$ and $Y(1)$, where $X^{\vee}$ and $Y$ are their cocharacter groups;
    \item $A$ and $A^*$ are two abelian $S$-schemes dual to each other;
    \item $v:X \longrightarrow A$ and $v^*:Y^{\vee} \longrightarrow A^*$ are two morphisms  of $S$-group schemes; and 
    \item $\psi$ is a trivialization of the pull-back $(v,v^*)^*{\mathcal{P}}_A$ via $(v,v^*)$ of the Poincar\' e biextension ${\mathcal{P}}_A$ of $(A,A^*)$.
\end{itemize}    

\pn To have the data $(X,Y^{\vee},A,A^*, v ,v^*,\psi)$ is equivalent to have the 1-motive $M=[X {\buildrel u \over \longrightarrow }G]$: In fact,
to have the semi-abelian $S$-scheme $G$ is the same thing as to have the morphism 
$v^*:Y^{\vee} \longrightarrow A^*$ (cf. ~\cite{SGA7} Expos\'e VIII 3.7, $G$ corresponds to the biextension $(id_A,v^*)^*{\mathcal{P}}_A$ of $(A,Y^{\vee})$ by ${\GG}_m$) and 
to have the morphism $u:X \longrightarrow G$ 
is equivalent to have the morphism $v:X \longrightarrow A$ and the trivialization $\psi$ of $(v,v^*)^*{\mathcal{P}}_A$ (the trivialization $\psi$ furnishes the lift $u:X \longrightarrow G$ of the morphism $v:X \longrightarrow A$).

\par Let $s: X \times Y^{\vee} \longrightarrow Y^{\vee} \times X$ be 
the morphism which permutates the factors. 
The 7-uplet $(Y^{\vee},X,A^*,A, v^* ,v,\psi \circ s)$ defines 
the so called \textbf{Cartier dual of $M$}: this is a 1-motive that we denote 
by $M^*$.

\begin{rem}
The pull-back $(v,v^*)^*{\mathcal{P}}_A$ by $(v,v^*)$
of the Poincar\' e biextension ${\mathcal{P}}_A$ of $(A,A^*)$ is a biextension
of $(X,Y^{\vee})$ by ${\GG}_m$. 
 According~\cite{SGA3} Expos\'e X Corollary 4.5, we can suppose 
that the character group $Y^{\vee}$ is constant, i.e. $\ZZ^{\rk Y^{\vee}}$ 
(if necessary we localize over $S$ for the \'etale topology).
Moreover since by~\cite{SGA7} Expos\'e VII (2.4.2) the category $\bBiext$ is additive in each variable, we have that
\[{\bBiext}(X,Y^{\vee};{\GG}_m) \cong {\bBiext}(X,{\ZZ}; Y(1)).\]
\pn We denote by 
\[((v,v^*)^*{\mathcal{P}}_A) \otimes Y\]
\pn the biextension of $(X,\ZZ)$ by $Y(1)$ corresponding to the 
 biextension $(v,v^*)^*{\mathcal{P}}_A$ through this equivalence of categories. 
The trivialization $\psi$ of $(v,v^*)^*{\mathcal{P}}_A$ defines a trivialization 
\[\psi \otimes Y \]
\pn  of $((v,v^*)^*{\mathcal{P}}_A) \otimes Y$, and vice versa.
\end{rem}

We can now give a more useful definition of a biextension of two 1-motives by a third one:

\begin{prop}\label{defbiext2}
Let $M_i=(X_i,Y_i^{\vee},A_i,A_i^*, v_i ,v_i^*,\psi_i)$ (for $i=1,2,3$) be a 1-motive. A biextension $B=(B,\Psi'_1, \Psi'_2,\Psi',\Lambda)$ of $(M_1,M_2)$ by $M_3$ consists of 
\begin{enumerate}
    \item a biextension of $B$ of $(A_1, A_2)$ by $Y_3(1);$
    \item a trivialization $\Psi'_1 $ (resp. $\Psi'_2$) of the biextension 
$(v_1, id_{A_2})^*B$ (resp.\\
 $(id_{A_1},v_2)^*B$) of $(X_1, A_2)$ by $Y_3(1)$ 
(resp. of $(A_1, X_2)$ by $Y_3(1)$) obtained as pull-back of the biextension $B$ via $(v_1,id_{A_2})$ (resp. via $(id_{A_1},v_2)$); 
    \item a trivialization $\Psi'$ of the biextension 
$(v_1, v_2)^*B$ of $(X_1, X_2)$ by $Y_3(1)$ obtained as pull-back of the biextension 
$B$ via $(v_1,v_2)$, which coincides with the trivializations 
induced by $\Psi'_1$ and $\Psi'_2$ over $X_1 \times X_2$, i.e. 
\[(v_1,id_{A_2})^*\Psi'_2=\Psi'=(id_{A_1},v_2)^*\Psi'_1; \]
    \item a morphism $\Lambda: (v_1, v_2)^*B \longrightarrow ((v_3,v_3^*)^*{\mathcal{P}}_{A_3}) \otimes Y_3 $ of trivial biextensions, with 
$\Lambda_{| Y_3(1)}$ equal to the the identity, such that the following diagram is commutative 

\begin{equation}\label{condbiext}
\begin{array}{ccc}
  Y_3(1) &= &  Y_3(1)\\
        \vert &     & \vert\\
(v_1, v_2)^*B &  \longrightarrow &((v_3,v_3^*)^*{\mathcal{P}}_{A_3}) \otimes Y_3 \\
{\scriptstyle \Psi'} \uparrow \downarrow~~ & &~~~~~~~~~\downarrow \uparrow {\scriptstyle \psi_3 \otimes Y_3}   \\
X_1 \times X_2 & \longrightarrow & X_3 \times \ZZ. \\
\end{array}
\end{equation}

\end{enumerate}    
\end{prop}

\begin{proof}
According to the main Theorem of~\cite{B2}, to have the biextension $B$ of $(A_1,A_2)$ by $Y_3(1)$ is equivalent to have the  biextension 
${\mathcal{B}}=\iota_{3\,*}(\pi_1,\pi_2)^*B$ of $(G_1, G_2)$ by $G_3$, where for $i=1,2,3,$ $\pi_i:G_i \longrightarrow A_i$ 
is the projection of $G_i$ over $A_i$ and $\iota_i:Y_i(1) \longrightarrow 
G_i$ is the inclusion of $Y_i(1)$ over $G_i.$ \\
The trivializations $(\Psi' ,\Psi'_1 ,\Psi'_2 )$ and $(\Psi,\Psi_1 ,\Psi_2 )$ determine each others.\\
By~\cite{B2}, both biextensions $(v_1, v_2)^*B$ and $((v_3,v_3^*)^*{\mathcal{P}}_{A_3}) \otimes Y_3$ are trivial.
Hence to have the morphism of $S$-group schemes $\lambda:X_1 \times X_2 \longrightarrow X_3$ is equivalent to have the morphism 
of trivial biextensions
$\Lambda: (v_1, v_2)^*B \longrightarrow ((v_3,v_3^*)^*{\mathcal{P}}_{A_3}) \otimes Y_3 $  with $\Lambda_{| Y_3(1)}$ equal to the identity. In particular, through this equivalence $\lambda$ corresponds to 
$\Lambda_{| X_1 \times X_2}$ and
to require that
$u_3 \circ \lambda: X_1 \times X_2 \longrightarrow G_3$ is compatible with the trivialization  $\Psi$ of $(u_1,u_2)^* {\mathcal{B}}$ corresponds to require the commutativity of the diagram (\ref{condbiext}). 
\end{proof}

\begin{rem}
The data (1), (2), (3) and (4) of definition \ref{defbiext} are equivalent respectively to the data (1), (2), (3) and (4) of Proposition \ref{defbiext2}.
\end{rem}

\subsection{Examples}\label{exambiext}
We conclude this chapter giving some examples of biextensions of 1-motives by 1-motives. We will use the more useful definition of biextensions furnishes by Proposition \ref{defbiext2}.

\textbf{(1)} Let $M=[0 \longrightarrow A]$ be an abelian $S$-scheme with Cartier dual
 $M^*=[0 \longrightarrow A^*]$ and let $W(1)$ be an $S$-torus.
A biextension of $(M,M^*)$ by $W(1)$ is 
\[(B,0,0,0,0)\]
\pn where $B$ is a biextension of $(A,A^*)$ by $W(1)$. In particular 
\textbf{the Poincar\'e biextension of 
$(M,M^*)$ by ${\ZZ}_(1)$} is the biextension
\[({\mathcal{P}}_A,0,0,0,0)\]
\pn where ${\mathcal{P}}_A$ is the Poincar\' e biextension of $(A,A^*)$.

\textbf{(2)} Let $M=(A,A^*, X,Y^{\vee}, v ,v^*,\psi)=[X \stackrel{u}{\longrightarrow }G]$
 be a 1-motive over $S$ and $M^*=[Y^{\vee} \stackrel{u^*}{ \longrightarrow }G^*]$ its Cartier dual.
If ${\mathcal{P}}_A$ denotes the Poincar\' e biextension of $(A,A^*)$, the semi-abelian $S$-scheme $G$
(resp. $G^*$) corresponds to the biextension $(id_A,v^*)^* {\mathcal{P}}_A$ of $(A,Y^{\vee})$ by ${\ZZ}(1)$
(resp. $(v,id_{A^*})^* {\mathcal{P}}_A$ of $(X,A^*)$ by ${\GG}_m$)
(cf.~\cite{SGA7} Expos\'e VIII 3.7). \textbf{The Poincar\'e biextension of 
$(M,M^*)$ by ${\ZZ}(1)$} is the biextension
\[({\mathcal{P}}_A, \psi_1, \psi_2,\psi,0)\]
\pn where $\psi_1$ is the trivialization of the biextension $(id_A,v^*)^* {\mathcal{ P}}_A$ which
defines the morphism $u: X \longrightarrow G$, and $\psi_2$ is the trivialization of the biextension $(v,id_{A^*})^* {\mathcal{P}}_A$ which defines the morphism $u^*: Y^{\vee} \longrightarrow G^*.$

\textbf{(3)} Let $M=[X \stackrel{u}{\longrightarrow } Y(1)]$ be the 1-motive without abelian part defined over $S$. Its Cartier dual is the 1-motive 
$M^*=[Y^{\vee} \stackrel{ u^*}{ \longrightarrow } X^{\vee}(1)],$ where $Y^{\vee}$ is the character group of the $S$-torus $Y(1)$ and $X^{\vee}(1)$ is the $S$-torus whose character group is $X$. \textbf{ The Poincar\'e biextension of $(M,M^*)$ by ${\ZZ}(1)$} is the biextension
\[({\mathbf{\underline 0}}, 0, 0,\psi,0)\]
\pn where ${\mathbf{\underline 0}}$ is the trivial biextension of $(0,0)$
by ${\ZZ}(1)$ and 
$\psi: X \times Y^{\vee} \longrightarrow {\ZZ} (1)$ is the biaddictive morphism defining $u:X \longrightarrow Y(1)$ and $u^*:Y^{\vee} \longrightarrow X^{\vee}(1)$. (We can view the trivialization $\psi$ as a biaddictive morphism from 
$ X \times Y^{\vee}$ to ${\ZZ} (1)$ because of the triviality of the biextension ${\mathbf{\underline 0}}$).

\textbf{(4)} Let $M=[X \stackrel{u}{\longrightarrow } Y(1)]$ and $[V \stackrel{t}{\longrightarrow } W(1)]$ be two 1-motives without abelian part 
defined over $S$, and let $M^*=[Y^{\vee} \stackrel{ u^*}{ \longrightarrow } X^{\vee}(1)]$ be the Cartier dual of $M$. A biextension of $(M,M^*)$ by $[V \stackrel{t}{\longrightarrow } W(1)]$ is
\[({\mathbf{\underline 0}},0,0,\Psi,\lambda)\]
\pn where ${\mathbf{\underline 0}}$ is the trivial biextension of $(0,0)$
by $W(1)$, $\Psi: X \times Y^{\vee} \longrightarrow W(1)$ is a biaddictive morphism and $\lambda:X \times Y^{\vee} \longrightarrow V$ a morphism of $S$-group schemes, 
such that the diagram

\begin{equation}\label{diagr:1}
\begin{array}{ccc}
  W(1) & = & W(1) \\
 \Psi \uparrow &  & \uparrow t \\
  X \times Y^{\vee} & \stackrel{\lambda}{\longrightarrow} & V \\
\end{array}
\end{equation}
 
\pn is commutative.(We can view the trivialization $\Psi$ as a
 biaddictive morphism from 
$ X \times Y^{\vee}$ to $W (1)$ because of the triviality of the biextension ${\mathbf{\underline 0}}$).

\textbf{(5)} Let $M=[X \stackrel{u}{\longrightarrow } Y(1)]$ be a 1-motive without abelian part defined over $S$, and $M^*=[Y^{\vee} \stackrel{ u^*}{ \longrightarrow } X^{\vee}(1)]$ its Cartier dual. Let $[V \stackrel{t}{\longrightarrow } A]$ be a 1-motive defined over $S$ without toric part.
A biextension of $(M,M^*)$ by $[V \stackrel{t}{\longrightarrow } A]$ is
\[({\mathbf{\underline 0}},0,0,0,\lambda)\]
\pn where ${\mathbf{\underline 0}}$ is the trivial biextension of $(0,0)$
by $0$ and $\lambda:X \times Y^{\vee} \longrightarrow V$ is a morphism of $S$-group schemes. Therefore 

\begin{equation}\label{biextMM[VA]}
{\bBiext}^1(M,M^*;[V \stackrel{t}{\longrightarrow } A]) = {\Hom}_{S-{\mathrm{schemes}}}(X \times Y^{\vee},V).
\end{equation}

\textbf{(6)} Let $M=[X \stackrel{u}{\longrightarrow } A \times Y(1)]=(X,Y^{\vee},A,A^*, v ,v^*=0,\psi)$ be a 1-motive defined over $S$, and $M^*=[Y^{\vee} \stackrel{ u^*}{ \longrightarrow } G^*]$ its Cartier dual. Let $[V \stackrel{t}{\longrightarrow } W(1)]$ be a 1-motive defined over $S$ without abelian part.
A biextension of $(M,M^*)$ by $[V \stackrel{t}{\longrightarrow } W(1)]$ is
\[(B,\Psi_1,\Psi_2,\Psi,\lambda)\]
\pn where 

\begin{itemize}
    \item $B$ is a biextension of $(A, A^*)$ by $W(1);$
    \item $\Psi_1 $ (resp. $\Psi_2$) is a trivialization of the biextension 
$(v, id_{A^*})^*B$
\pn (resp. $(id_{A},v^*)^*B$) of $(X, A^*)$ by $W(1)$ 
(resp. of $(A, Y^{\vee})$ by $W(1)$) obtained as pull-back of the biextension $B$ via $(v,id_{A^*})$ (resp. via $(id_{A},v^*)$). Since $v^*=0$, the biextension 
$(id_{A},v^*)^*B$ of $(A, Y^{\vee})$ by $W(1)$ is trivial; 
    \item $\Psi$ is a trivialization of the biextension 
$(v, v^*)^*B$ of $(X, Y^{\vee})$ by $W(1)$ obtained as pull-back of the biextension 
$B$ via $(v,v^*)$, which coincides with the trivializations 
induced by $\Psi_1$ and $\Psi_2$ over $X \times Y^{\vee}$. Also the biextension
$(v, v^*)^*B$ of $(X, Y^{\vee})$ by $W(1)$ is trivial and hence we can view 
$\Psi$ as a biaddictive morphism from $X \times Y^{\vee}$ to $W(1)$;
    \item $\lambda: X \times Y^{\vee} \longrightarrow V $ is a morphism of $S$-group schemes such that the following diagram is commutative 

\begin{equation}\label{diagr:2}
\begin{array}{ccc}
  W(1) &= &  W(1)\\
{\scriptstyle \Psi} \uparrow ~~ & & \uparrow {\scriptstyle t}   \\
X \times Y^{\vee} & \longrightarrow & V.\\
\end{array}
\end{equation}

\end{itemize}

\textbf{(7)} Let $M=[X \stackrel{u}{\longrightarrow } A \times Y(1)]=(X,Y^{\vee},A,A^*, v ,v^*=0,\psi)$ be a 1-motive defined over $S$, and $M^*=[Y^{\vee} \stackrel{ u^*}{ \longrightarrow } G^*]$ its Cartier dual. Let $[V \stackrel{t}{\longrightarrow } A']=(V,0,A',{A'}^*,t,0,\psi')$ be a 1-motive defined over $S$ without toric part.
 A biextension of $(M,M^*)$ by $[V \stackrel{t}{\longrightarrow } A']$ is
\[({\mathbf{\underline 0}},0,0,\Psi,\lambda)\]
\pn where ${\mathbf{\underline 0}}$ is the trivial biextension of $(A,A^*)$
by $0$, $\Psi: X \times Y^{\vee} \longrightarrow X \times Y^{\vee}$ is a biaddictive morphism and $\lambda:X \times Y^{\vee} \longrightarrow V$ a morphism of $S$-group schemes, 
such that the diagram

\begin{equation}\label{diagr:3}
\begin{array}{ccc}
  X \times Y^{\vee} & \stackrel{\lambda}{\longrightarrow} &  V \times {\ZZ}\\
 \Psi \uparrow ~~~~ &  & ~~~~~~~ \uparrow \psi'\otimes 0 \\
  X \times Y^{\vee} & \stackrel{\lambda}{\longrightarrow} & V \times {\ZZ}\\
\end{array}
\end{equation}
 
\pn is commutative.(We can view the trivialization $\Psi$ as a
 biaddictive morphism from 
$ X \times Y^{\vee}$ to $X \times Y^{\vee}$ because of the triviality of the biextension ${\mathbf{\underline 0}}$ and for the meaning of $\psi' \otimes 0$ see (\ref{condbiext})). 

%----------------------------------------------------------------------------------
\section{Some morphisms of 1-motives}
In this chapter we suppose that 1-motives 
over $S$ generate a Tannakian category over a field $K$. 

Let ${\mathcal{M}}(S)$ be the Tannakian category generated by 
1-motives over $S$. 
The unit object of ${\mathcal{M}}(S)$ is the 1-motive 
${\ZZ}(0)=[{\ZZ}  \longrightarrow 0]$.
We denote by $M^{\vee} \cong {\underline {\Hom}}(M,{\ZZ}(0)) $
the dual of the 1-motive $M$ and
 $ev_M: M \otimes M^{\vee} \longrightarrow {\ZZ}(0)$ and
$\delta_M: {\ZZ}(0) \longrightarrow M^{\vee} \otimes M$ the morphisms of 
 ${\mathcal{M}}(S)$ which characterize  $M^{\vee}$ (cf.~\cite{D2} (2.1.2)).
 The Cartier dual of $M$ is the 1-motive

\begin{equation}\label{du}
M^*=M^{\vee} \otimes {\ZZ}(1)
\end{equation}

\pn where ${\ZZ}(1)$ is the $S$-torus with cocharacter group $\ZZ$.

\subsection{Motivic remarks about morphisms of 1-motives}

\begin{lem}\label{georem0}
Let $X_i$ (for $i=1,2,3$) be a $S$-group scheme which is locally 
for the \'etale topology a constant group scheme defined by a
 finitely generated free $\ZZ$-module, and let $G_i$ (for $i=1,2,3$) be a semi-abelian $S$-scheme.
In the category ${\mathcal{M}}(S),$ there are no morphisms from the tensor product 
$X_1 \otimes X_2$ to $G$, i.e.,
\[ {\Hom}_{{\mathcal{M}}(S)}(X_1 \otimes X_2,G_3) =0\]
\pn and there are no morphisms from the tensor product $G_1 \otimes G_2$ to $X_3$, i.e.,
\[ {\Hom}_{{\mathcal{M}}(S)}(G_1 \otimes G_2,X_3) =0.\]
\end{lem}

\begin{proof} This lemma follows from the fact that morphisms of motives have to respect weights. In fact the pure motive $X_1 \otimes X_2$ has weight 0 and
 the mixed motive $G$ has weight smaller or equal to - 1, and the mixed motive 
 $G_1 \otimes G_2$ has weight smaller or equal to -2 and the pure motive $X_3$
has weight 0. 
\end{proof}

\begin{prop}\label{georem1}
Let $G_i$ (for $i=1,2,3$) be a semi-abelian $S$-scheme, i.e.
an extension of an abelian $S$-scheme $A_i$ by a $S$-torus $Y_i(1)$. 
We have that:

\begin{eqnarray}
  {\Hom}_{{\mathcal{M}}(S)}(G_1 \otimes G_2,Y_3(1)) &=& {\Hom}_{{\mathcal{M}}(S)}(A_1 \otimes A_2,Y_3(1)) \\
  {\Hom}_{{\mathcal{M}}(S)}(G_1 \otimes G_2,G_3) &=& {\Hom}_{{\mathcal{M}}(S)}(G_1 \otimes G_2,Y_3(1)) 
\end{eqnarray}  
\end{prop}

\begin{proof} The proof of these equalities is based on the fact that morphisms of motives have to respect weights.\\
For $i=1,2,3,$ denote by $\pi_i:G_i \longrightarrow 
A_i$ the projection of $G_i$ over $A_i$ and $\iota_i:Y_i(1) \longrightarrow 
G_i$ the inclusion of $Y_i(1)$ in $G_i.$
\pn Thanks to the projection $(\pi_1, \pi_2),$ each morphism 
$A_1 \otimes A_2 \longrightarrow Y_3(1)$ can be lifted to a morphism
from $G_1 \otimes G_2$ to $Y_3(1)$. In the other hand, since the motive 
$Y_1(1) \otimes Y_2(1)$ has weight -4, each morphism 
$G_1 \otimes G_2 \longrightarrow Y_3(1)$ factorizes through
$A_1 \otimes A_2$.\\
Through the inclusion $\iota_3$, each morphism 
$G_1 \otimes G_2 \longrightarrow Y_3(1)$ defines a morphism
$G_1 \otimes G_2 \longrightarrow G_3$. In the other hand, since the motive 
$G_1 \otimes G_2$ has weight less or equal to -2, each morphism 
$G_1 \otimes G_2 \longrightarrow G_3$ factorizes through $Y_3(1)$.
\end{proof}

Therefore we have:

\begin{cor}\label{georem2}
In ${\mathcal{M}}(S)$, a morphism from the tensor product of two semi-abelian 
$S$-schemes to a semi-abelian $S$-scheme is a morphism from the tensor product 
of the underlying abelian $S$-schemes to the underlying $S$-torus:
\[{\Hom}_{{\mathcal{M}}(S)}(G_1 \otimes G_2,G_3) ={\Hom}_{{\mathcal{M}}(S)}(A_1 \otimes A_2,Y_3(1)). \]
\end{cor}

\subsection{Morphisms from a finite tensor products of 1-motives to a 1-motive}

\begin{defn}\label{defmor}
In the category ${\mathcal{M}}(S),$ \textbf{the morphism 
$M_1 \otimes M_2 \longrightarrow M_3$ from the tensor product of two 1-motives
to a third 1-motive} is an isomorphism class of biextensions of $(M_1, M_2)$ 
by $M_3$ (cf. (\ref{defbiext})). We define
\begin{equation}
{\Hom}_{{\mathcal{M}}(S)}(M_1 \otimes M_2,M_3) = {\Biext}^1(M_1,M_2;M_3) 
\end{equation}
\end{defn}

In other words, \emph{the biextensions of two 1-motives by a 1-motive are the ``geometrical interpretation'' of the morphisms of ${\mathcal{M}}(S)$ from the tensor product of two 1-motives to a 1-motive.}
\vskip 0.5 true cm

\begin{rem} The set ${\Hom}_{{\mathcal{M}}(S)}(M_1 \otimes M_2,M_3)$ is a group.
\end{rem}

\begin{rem}
The confrontation of the main Theorem of~\cite{D2} with
Corollary \ref{georem2}, shows that this definition is compatible with what we know about morphisms between 1-motives.
\end{rem}

\begin{defn} Let $M_i$ and $M'_i$ (for $i=1,2,3$) be 1-motives over $S$.
The notion \ref{defmorphbiext} of morphisms of biextensions defines 
a \textbf{morphism from the group of morphisms of ${\mathcal{M}}(S)$ 
from $ M_1 \otimes M_2$ to $ M_3,$
to the group of morphisms of ${\mathcal{M}}(S)$ from 
$ M'_1 \otimes M'_2 $ to $ M'_3$}, i.e.
\[ {\Hom}_{{\mathcal {M}}(S)}(M_1 \otimes M_2,M_3) \longrightarrow {\Hom}_{{\mathcal{M}}(S)}(M'_1 \otimes M'_2,M'_3).\]
\end{defn}

\begin{rem}
Using the notaions of \ref{defmorphbiext}, if we denote $b$ the morphism
$ M_1 \otimes M_2 \longrightarrow M_3$ corresponding to the biextension 
$({\mathcal{B}},\Psi_{1}, \Psi_{2},\lambda)$ of $(M_1,M_2)$ by $M_3$
and by $b'$ the morphism
$ M'_1 \otimes M'_2 \longrightarrow  M'_3$ corresponding to the biextension 
$({\mathcal{B}}',\Psi'_{1}, \Psi'_{2},\lambda')$ of $(M'_1,M'_2)$ by $M'_3$,
the morphism 

\[(F,\Upsilon_1,\Upsilon_2,\Upsilon,g_3):
({\mathcal{B}},\Psi_{1}, \Psi_{2},\lambda)
\longrightarrow ({\mathcal{B}}',\Psi'_{1}, \Psi'_{2},\lambda')\]

\pn of biextensions defines the vertical arrows of the following diagram of morphisms of 
${\mathcal{M}}(S)$

\[\begin{array}{ccc}
  M_1 \otimes M_2  & \stackrel{b}{\longrightarrow} &M_3  \\
 \downarrow  &  &  \downarrow \\
  M'_1 \otimes M'_2 & \stackrel{b'}{\longrightarrow} & M'_3. \\
\end{array}\]

\pn It is clear now
why from the data $(F,\Upsilon_1,\Upsilon_2,\Upsilon,g_3)$ we get a
morphism of ${\mathcal {M}}(S)$ from $M_3$ to $M'_3$ as
remarked in \ref{remmorphbiext}. Moreover
since $ M_1 \otimes {\ZZ}(0)$, $ M'_1 \otimes {\ZZ}(0)$,
$ {\ZZ}(0) \otimes  M_2$ and $ {\ZZ}(0) \otimes M'_2$, are sub-1-motives 
of the motives $M_1 \otimes M_2$ and $M'_1 \otimes M'_2$, it is clear that
from the data 
$(F,\Upsilon_1,\Upsilon_2,\Upsilon,g_3)$ we get morphisms 
 from $M_1$ to $M'_1$ and from $M_2$ to $M'_2$ as
remarked in \ref{remmorphbiext}.
\end{rem}

\begin{lem}\label{otimesM}
Let $l$ and $i$ be positive integers and let $M_j=[X_j \stackrel{u_j}{\longrightarrow} G_j]$ (for $j=1, \dots,l$) be a 1-motive defined over $S$. Denote by $M_0$ or $X_0$ the 1-motive ${\ZZ}(0)=[{\ZZ} \longrightarrow 0]$. If $i \geq 1$ and $l +1 \geq i$, the motive $ \otimes^l_{j=1}M_j /{\W}_{-i}(\otimes^l_{j=1}M_j)$ is isogeneous to the motive 

\begin{equation}\label{formulaotimes}
  \sum
 \Big(\otimes_{k \in \{\nu_0, \dots,\nu_{l-i+1}\}}X_k \Big) \bigotimes 
 \Big(\otimes_{j\in \{\iota_0, \dots, \iota_{i-1}\}} M_j /{\W}_{-i}(\otimes_{j\in \{\iota_0 , \dots, \iota_{i-1}\}}M_j) \Big)
\end{equation}

\pn where the sum is taken over all the $(l-i+1)$-uplets $\{\nu_0, \dots,\nu_{l-i+1}\}$ and all the $(i-1)$-uplets $\{\iota_0, \dots, \iota_{i-1}\}$
of $\{0,1, \cdots,l\}$ such that $ \{\nu_0, \dots,\nu_{l-i+1}\}\cap \{\iota_0, \dots, \iota_{i-1}\} = \emptyset$, 
$\nu_0 < \nu_1 < \dots < \nu_{l-i+1}$, $ \iota_0 < \iota_1 < \dots < \iota_{i-1},$ 
$\nu_a \not= \nu_b$ and $\iota_c \not= \iota_d$, for all $a,b \in \{0, \dots, l-i+1 \}, a \not=b$ and $c,d \in \{0, \dots, i-1 \}, c \not=d$.
\end{lem}

\begin{proof} 
1-motives $M_j$ are composed by pure motives of weight 0 (the lattice part $X_j$), 
-1 (the abelian part $A_j$) and -2 (the toric part $Y_j(1)$).
Consider the pure motive ${\Gr}^{\W}_{-i} (\otimes^l_{j=1}M_j)$: it is a finite sum of tensor products of $l$ factors of weight 0, -1 other -2.
If $i=l$ the tensor product 
\[A_1 \otimes A_2 \otimes \dots \otimes A_l\]
\pn contains no factors of weight 0. For each $i$ strictly bigger than $l$, it is also easy to construct a tensor product of $l$ factors whose total weight is $-i$ and in which no factor has weight 0 (for example if 
$i=l+2$ we take
\[Y_1(1) \otimes Y_2(1) \otimes A_3 \otimes \dots \otimes A_l.)\]
\pn However if $i$ is strictly smaller than $l$, 
in each of these tensor products of $l$ factors, there is at least one factor of weight 0, i.e. one of the $X_j$ for $j=1, \dots,l$.\\
Now fix a $i$ strictly smaller than $l$. The tensor products where there are less factors of weight 0 are exactly those where there are more factors of weight -1. Hence in the pure motive ${\Gr}_{-i} (\otimes^l_{j=1}M_j)$, the tensor products with less factors of weight 0 are of the type
\[ X_{\nu_1} \otimes \dots \otimes X_{\nu_{l-i}} \otimes A_{\iota_1} \otimes \dots \otimes A_{\iota_i}\]
\pn Thanks to these observations, the conclusion is clear.\\
Remark that we have only an isogeny because in the 1-motive (\ref{formulaotimes})
the factor

\[X_{\nu_1} \otimes X_{\nu_2} \otimes \dots \otimes X_{\nu_p} \otimes  
{\mathcal{Y}}_{\iota_1} \otimes {\mathcal{Y}}_{\iota_2} \otimes  \dots \otimes 
{\mathcal{Y}}_{\iota_{l-p}}\]

\pn appears with multiplicity ``$p+m$'' where $m$ is the number of 
${\mathcal{Y}}_{\iota_q}$ (for $q=1, \dots,l-p$) which are of weight 0, instead of appearing only once like in the 1-motive \\
$ \otimes_{j}M_j /{\W}_{-i}(\otimes_{j}M_j)$.
In particular for each $i$ we have that

\begin{eqnarray}
\nonumber  {\Gr}^{\W}_0 \Big(
\sum (\otimes_{k }X_k ) \otimes (\otimes_{j} M_j /{\W}_{-i}) \Big) &=& l ~{\Gr}^{\W}_0 \Big( \otimes_{j}M_j /{\W}_{-i} \Big) \\
\nonumber  {\Gr}^{\W}_{-1} \Big(
\sum (\otimes_{k }X_k ) \otimes (\otimes_{j} M_j /{\W}_{-i}) \Big) &=& (l-1) ~{\Gr}^{\W}_{-1} \Big( \otimes_{j}M_j /{\W}_{-i} \Big) \end{eqnarray}  
\end{proof}

We will denote by ${\mathcal{M}}^{\mathrm{iso}}(S)$ the Tannakian category generated by the iso-1-motives, i.e. by 1-motives modulo isogenies.

\begin{thm} \label{thmotimes}
Let $M$ and $M_1, \dots, M_l$ be 1-motives over $S$.
In the category ${\mathcal{M}}^{\mathrm{iso}}(S),$ the morphism 
$\otimes^l_{j=1}M_j \longrightarrow M$ from a finite tensor product of 1-motives
to a 1-motive is the sum of copies of isomorphism classes of biextensions of $(M_i, M_j)$ by $M$ for $i,j=1, \dots l$ and $i \not= j$. 
We have that
\[{\Hom}_{{\mathcal{M}}^{\mathrm{iso}}(S)}(\otimes^l_{j=1}M_j,M) =\sum_{i,j \in \{1, \dots,l\} \atop i \not= j} {\Biext}^1(M_i,M_j;M) \]
\end{thm}

\begin{proof} Because morphisms between motives have to respect weights,
the non trivial components of the morphism $\otimes^l_{j=1}M_j \longrightarrow M$ 
are the one of the morphism
\[ \otimes^l_{j=1}M_j \Big/ {\W}_{-3}(\otimes^l_{j=1}M_j) \longrightarrow M.\]
\pn Using the equality obtained in Lemma \ref{otimesM} with $i=-3$, we can write explicitly this last morphism in the following way
\[ \sum_{ \iota_1 < \iota_2 ~\mathrm{and}~ \nu_1< \dots < \nu_{l-2}
 \atop \iota_1 , \iota_2 \notin \{\nu_1, \dots,\nu_{l-2}\} }
 X_{\nu_1}\otimes \dots \otimes X_{\nu_{l-2}}\otimes ( M_{\iota_1} \otimes M_{\iota_2} /{\W}_{-3}( M_{\iota_1} \otimes M_{\iota_2})) \longrightarrow M.\]
Since ``to tensorize a motive by a motive of weight 0'' means to take a certain number of copies of the motive, from definition \ref{defmor} we get the expected conclusion.
\end{proof}

\subsection{Examples}
Now we give some examples of morphisms from the tensor product of two 1-motives 
to a 1-motive. According to our definition \ref{defmor} of morphisms of ${\mathcal{M}}(S)$, these examples are in parallel with the examples \ref{exambiext}.

\textbf{(1)} Let $M=[0 \longrightarrow A]$ be an abelian $S$-scheme with Cartier dual
 $M^*=[0 \longrightarrow A^*]$ and let $W(1)$ be an $S$-torus.
By definition, the biextension $(B,0,0,0,0)$ of $(M,M^*)$ by $W(1)$ is 
a morphism of ${\mathcal{M}}(S)$ from $A \otimes A^*$ to $W(1)$.
In particular, the Poincar\'e biextension ${\mathcal{P}}_A$ 
of $(A,A^*)$ by ${\ZZ}(1)$ is the classical \textbf{Weil pairing of $A$}
\[{\mathcal{P}}_A: A \otimes A^* \longrightarrow {\ZZ}(1).\]

\textbf{(2)} Let $M$ be a 1-motive over $S$ and $M^*$ its Cartier dual. 
The Poincar\'e biextension ${\mathcal{P}}_M$ 
of $(M,M^*)$ by ${\ZZ}(1)$ (cf. \ref{exambiext})
 is the so called \textbf{Weil pairing of $M$}
\[{\mathcal{P}}_M:M \otimes M^* \longrightarrow {\ZZ}(1),\]
\pn which expresses the Cartier duality between $M$ and $M^*$. 

The valuation map $ev_M: M \otimes M^{\vee} \longrightarrow {\ZZ}(0)$
of $M$ expresses the duality between $M$ and $M^{\vee}$ as objects of the Tannakian category  ${\mathcal{M}}(S)$. By (\ref{du})
the evaluation map $ev_M$ is the twist by ${\ZZ}(-1)$ of the Weil pairing of $M$:

\[ev_M ={\mathcal{P}}_M \otimes {\ZZ}(-1): M \otimes M^{\vee} \longrightarrow {\ZZ}(0).\]

 \textbf{(3)} Let $M$ be the 1-motive without abelian part $[X \stackrel{u}{ \longrightarrow} Y(1)]$ defined over S. Its Cartier dual is the 1-motive 
$M^*=[Y^{\vee} \stackrel{u^*} {\longrightarrow} X^{\vee}(1)].$
We will compute the morphism from 
the tensor product $M \otimes M^* $ to the torus $ {\ZZ}(1)$, which is by definition 
the isomorphism class of the Poincar\'e biextension of $(M,M^*)$ by ${\ZZ}(1).$
Regarding $M$ and $M^*$ as complexes of commutative groups $S$-schemes concentrated 
in degree 0 and -1, their tensor product is the complex 

\begin{equation}\label{MotimesMdu1}
[X \otimes Y^{\vee} 
\stackrel{\scriptscriptstyle{(-id_{\scriptscriptstyle X} \otimes u^*,u \otimes id_{\scriptscriptstyle Y^{\vee}})}}{\longrightarrow} 
X \otimes X^{\vee}(1) + Y(1) \otimes Y^{\vee} \stackrel{\scriptscriptstyle{u \otimes id_{X^{\vee}}+id_{Y(1)}\otimes u^*}}{\longrightarrow} Y(1)\otimes X^{\vee}(1)]
\end{equation}

\pn where $X \otimes Y^{\vee}$ is a pure motive of weight 0, $X \otimes X^{\vee}(1) + Y(1) \otimes Y^{\vee}$ is a pure motive of weight -2 and $Y(1)\otimes X^{\vee}(1)$
is a pure motive of weight -4.
Because of the weights, the only non-trivial component of a morphism from 
the tensor product $M \otimes M^* $ to the torus $ {\ZZ}(1)$ is 

\begin{equation}\label{G-m}
X \otimes X^{\vee}(1) + Y(1) \otimes Y^{\vee} \longrightarrow {\ZZ}(1) 
\end{equation}

\pn By example \ref{exambiext} (3), the Poincar\'e biextension of $(M,M^*)$ by ${\ZZ}(1)$ is the biextension $(0, 0, 0,\psi,0)$,
where $\psi: X \times Y^{\vee} \longrightarrow {\ZZ} (1)$ is the biaddictive morphism defining $u:X \longrightarrow Y(1)$ and $u^*:Y^{\vee} \longrightarrow X^{\vee}(1)$. 
This biaddictive morphism $\psi$ defines the only non-trivial component (\ref{G-m})
of $M \otimes M^* \longrightarrow {\ZZ}(1)$ through the following commutative diagram

\[\begin{array}{ccc}
  X \otimes Y^{\vee} &  &  \\
 {\scriptscriptstyle{(-id_{\scriptscriptstyle X} \otimes u^*,u \otimes id_{\scriptscriptstyle Y^{\vee}})}} \downarrow ~~~~~~~~~~~~& \searrow^\psi &  \\
 X \otimes  X^{\vee}(1) + Y(1) \otimes Y^{\vee} & \longrightarrow & {\ZZ}(1)  \\
\end{array}\]

\textbf{(4)} Let $M=[X \stackrel{u}{\longrightarrow } Y(1)]$ and $[V \stackrel{t}{\longrightarrow } W(1)]$ be two 1-motives without abelian part 
defined over $S$, and let $M^*=[Y^{\vee} \stackrel{ u^*}{ \longrightarrow } X^{\vee}(1)]$ be the Cartier dual of $M$.
According (\ref{MotimesMdu1}), the only non trivial components of a morphism from 
the tensor product $M \otimes M^* $ to $[V \stackrel{t}{\longrightarrow } W(1)]$
are 

\begin{eqnarray}
\label{poids0:1} X \otimes Y^{\vee} &\longrightarrow&  V \\
\label{poids-2:1} X \otimes X^{\vee}(1) + Y(1) \otimes Y^{\vee} &\longrightarrow& W(1) 
\end{eqnarray}  

\pn By example \ref{exambiext} (4), a biextension of $(M,M^*)$ by $[V \stackrel{t}{\longrightarrow } W(1)]$ is $({\mathbf{\underline 0}},0,0,\Psi,\lambda)$
where ${\mathbf{\underline 0}}$ is the trivial biextension of $(0,0)$
by $W(1)$, $\Psi: X \times Y^{\vee} \longrightarrow W(1)$ is a biaddictive morphism and $\lambda:X \times Y^{\vee} \longrightarrow V$ a morphism of $S$-group schemes, 
such that the diagram (\ref{diagr:1}) is commutative. The morphism $\lambda$ defines the non-trivial
component (\ref{poids0:1}) between motives of weight 0.
Through the commutative diagram

\[\begin{array}{ccc}
  X \otimes Y^{\vee} &  &  \\
 {\scriptscriptstyle{(-id_{\scriptscriptstyle X} \otimes u^*,u \otimes id_{\scriptscriptstyle Y^{\vee}})}} \downarrow & \searrow^\Psi &  \\
 X \otimes  X^{\vee}(1) + Y(1) \otimes Y^{\vee} & \longrightarrow & W(1)  \\
\end{array}\]

\pn the morphism $\Psi$ defines the non-trivial
component (\ref{poids-2:1}) between motives of weight -2. The commutativity 
of the diagram (\ref{diagr:1}) impose the compatibility between these two non trivial components.

\textbf{(5)} Let $M=[X \stackrel{u}{\longrightarrow } Y(1)]$ be a 1-motive without abelian part defined over $S$, and $M^*=[Y^{\vee} \stackrel{ u^*}{ \longrightarrow } X^{\vee}(1)]$ its Cartier dual. Let $[V \stackrel{t}{\longrightarrow } A]$ a 1-motive defined over $S$ without toric part.
Because of the weights, according (\ref{MotimesMdu1}) the only non trivial component of a morphism from 
the tensor product $M \otimes M^* $ to $[V \stackrel{t}{\longrightarrow } A]$
is 

\begin{equation}\label{poids0:2} 
X \otimes Y^{\vee} \longrightarrow  V 
\end{equation}  

\pn By example \ref{exambiext} (4), a biextension of $(M,M^*)$ by $[V \stackrel{t}{\longrightarrow } A]$ is $({\mathbf{\underline 0}},0,0,\lambda)$
where ${\mathbf{\underline 0}}$ is the trivial biextension of $(0,0)$
by $0$ and $\lambda:X \times Y^{\vee} \longrightarrow V$ is a morphism of $S$-group schemes. This morphism $\lambda$ defines the only non trivial component 
(\ref{poids0:2}) of the morphism $M \otimes M^* \longrightarrow [V \stackrel{t}{\longrightarrow } A]$. In particular by (\ref{biextMM[VA]})
\[ {\Hom}_{\mathcal{M}(S)}(M\otimes M^*,[V \stackrel{t}{\longrightarrow } A]) = {\Hom}_{S-{\mathrm{schemes}}}(X \times Y^{\vee},V).\]

\textbf{(6)}  Let $M=[X \stackrel{u}{\longrightarrow } A \times Y(1)]=(X,Y^{\vee},A,A^*, v ,v^*=0,\psi)$ be a 1-motive defined over $S$, and $M^*=[Y^{\vee} \stackrel{ u^*}{ \longrightarrow } G^*]$ its Cartier dual. Let $[V \stackrel{t}{\longrightarrow } W(1)]$ be a 1-motive defined over $S$ without abelian part. Because of the weights, the only non trivial components of a morphism from 
the tensor product $M \otimes M^* $ to $[V \stackrel{t}{\longrightarrow } W(1)]$
are 

\begin{eqnarray}
\label{poids0:3} X \otimes Y^{\vee} &\longrightarrow&  V \\
\label{poids-2:3} A \otimes A^* + X \otimes X^{\vee}(1) + Y(1) \otimes Y^{\vee} &\longrightarrow& W(1) 
\end{eqnarray}  

\pn By example \ref{exambiext} (6), a biextension of $(M,M^*)$ by $[V \stackrel{t}{\longrightarrow } W(1)]$ is $(B,\Psi_1,\Psi_2,\Psi,\lambda)$
 where $B$ is a biextension of $(A, A^*)$ by $W(1),$
$\Psi:X \times Y^{\vee} \longrightarrow W(1) $ is a biaddictive morphism and 
$\lambda: X \times Y^{\vee} \longrightarrow V $ is a morphism of $S$-group schemes such that the diagram (\ref{diagr:2}) is commutative. The morphism $\lambda$ defines the non-trivial
component (\ref{poids0:3}) between motives of weight 0.
Through the commutative diagram

\[\begin{array}{ccc}
  X \otimes Y^{\vee} &  &  \\
  \downarrow & \searrow^\Psi &  \\
 A \otimes A^* +X \otimes  X^{\vee}(1) + Y(1) \otimes Y^{\vee} & \longrightarrow & W(1)  \\
\end{array}\]

\pn the morphism $\Psi$ and the biextension $B$ define the non-trivial
component (\ref{poids-2:3}) between motives of weight -2: to be more precise, the biextension $B$ is the component $A \otimes A^* \longrightarrow W(1)$ and the 
morphism $\Psi$ determines the component $ X \otimes  X^{\vee}(1) + Y(1) \otimes Y^{\vee}  \longrightarrow  W(1).$
The commutativity of the diagram (\ref{diagr:2}) imposes the compatibility between these two non trivial components.
 Remark that the factor $A \otimes A^*$ plays no role in the diagram (\ref{diagr:2})
because $v^*=0$.

\textbf{(7)} Let $M=[X \stackrel{u}{\longrightarrow } A \times Y(1)]=(X,Y^{\vee},A,A^*, v ,v^*=0,\psi)$ be a 1-motive defined over $S$, and $M^*=[Y^{\vee} \stackrel{ u^*}{ \longrightarrow } G^*]$ its Cartier dual. Let $[V \stackrel{t}{\longrightarrow } A']=(V,0,A',{A'}^*,t,0,\psi')$ be a 1-motive defined over $S$ without toric part.
 Because of the weights, the only non trivial components of a morphism from 
the tensor product $M \otimes M^* $ to $[V \stackrel{t}{\longrightarrow } A']$
are 

\begin{eqnarray}
\label{poids0:4} X \otimes Y^{\vee} &\longrightarrow&  V \\
\label{poids-1:4} X \otimes A^* + A \otimes Y^{\vee}(1) &\longrightarrow& A' 
\end{eqnarray}  

\pn By example \ref{exambiext} (7), a biextension of $(M,M^*)$ by $[V \stackrel{t}{\longrightarrow } A]$ is $({\mathbf{\underline 0}},0,0,\Psi,\lambda)$
where ${\mathbf{\underline 0}}$ is the trivial biextension of $(A,A^*)$
by $0$, $\Psi: X \times Y^{\vee} \longrightarrow X \times Y^{\vee}$ is a biaddictive morphism and $\lambda:X \times Y^{\vee} \longrightarrow V$ a morphism of $S$-group schemes, such that the diagram (\ref{diagr:3}) is commutative.
The morphism $\lambda$ defines the component (\ref{poids0:4}) between motives of weight 0 and the biaddictive morphism 
$\Psi$ determines the component (\ref{poids-1:4}). 
The commutativity of the diagram (\ref{diagr:3}) imposes the compatibility between these two non trivial components.

% ------------------------------------------------------------------------

% ---------------------------------------------------------------------------
\end{document}